\documentclass[10pt]{article}
\usepackage{amstext}
\usepackage[fleqn]{amsmath}
\usepackage{amsfonts}
\usepackage{amssymb}
\usepackage{amsthm}
\usepackage{newlfont}
\usepackage{graphicx}
\usepackage{sectsty}
\theoremstyle{plain} \theoremstyle{theorem}
\newtheorem{thm}{Theorem}[section]
\theoremstyle{example}

\theoremstyle{corollary}
\newtheorem{cor}{Corollary}[section]
\theoremstyle{lemma}
\newtheorem{lemma}{Lemma}[section]
\theoremstyle{proposition}

\theoremstyle{axiom}

\theoremstyle{notation}

\theoremstyle{fact}

\theoremstyle{definition}
\newtheorem{defn}{Definition}[section]
\theoremstyle{remark}

\numberwithin{equation}{section}
\begin{document}
\title{On $q$-Bessel matrix polynomials}
\author{Ayman Shehata \thanks{%
E-mail: aymanshehata@science.aun.edu.eg}
,\; M. Tawfik \thanks{%
E-mail: m.tawfek@aun.edu.eg}
, Ayman M. Mahmoud \thanks{%
E-mail: ayman27@sci.nvu.edu.eg}\; and
Nada Mostafa \thanks{%
E-mail: nadamoustafa@sci.nvu.edu.eg}\\
{\small $^{\ast}$Department of Mathematics, Faculty of Science, Assiut University, Assiut 71516, Egypt.}\\
{\small $^{\dag,\ddag,\S}$ Department of Mathematics and computer, Faculty of Science,}\\
{\small New Valley University, El-khargah 72511, Egypt.}}
\date{}
\maketitle{}

\begin{abstract}
The aim of the present study is to establish some properties for $q$-Bessel matrix polynomials such as several $q$-differential matrix equation, $q$-differential matrix relations and $q$-recurrence matrix relations, and integral representation, $q$-Laplace and $q$-Mellin transforms with the help of $q$-Analysis. Furthermore, we give connections between $q$-Horn's matrix functions of two variables and $q$-Bessel matrix polynomials are given.
\end{abstract}
\section{Introduction}
Most of the special functions (SFs) are shown as a solution to ordinary or integral equations, so they perform extremely importantly in the investigation of SF as well as in the formalism of mathematical physics. In recent years the matrix version of some fundamental (or $q$-)SFs have been studied, mainly because of their applications to physics, statistics, wave mechanics, quantum mechanics, engineering and lie groups and Lie theory. Abdi \cite{ab}, Purohitands and Kalla \cite{puk} studied $q$-Laplace transforms of $q$-Bessel functions. Abdi, Exton and Srivastava investigated the $q$‐Analogue of Bessel polynomials in \cite{ab, exs}. In \cite{be}, Belgacem demonstrated the Sumudu transform of Bessel functions and applications. In \cite{uc}, Ucar determined the $q$‐Sumudu transforms of the $q$‐Analogues of Bessel functions. In \cite{ch}, Chen and Feng discussed the group theoretic origins of certain generating functions of Bessel polynomials. In \cite{dw1, dw2}, Dwivedi development the special matrix functions, hypergeometric function, $q$-Special functions and evaluated matrix $q$-Kummer equation and its solutions. In \cite{dwsh1}, Dwivedi and Sahai have obtained the matrix analogue of $q$-hypergeometric and $q$-Appell functions of two variables. In \cite{dwsh2, dwsh3}, Dwivedi and Sahai studied of the $q$-hypergeometric matrix functions of two and several variables and their $q$-difference matrix  equations. In \cite{dwsn}, Dwivedi and Sanjhira investigated the $q$-Hypergeometric matrix function and its $q$-Fractional calculus. Exton discussed a $q$-Analogue of Bessel-Clifford equation and $q$-Hypergeometric functions and applications in \cite{ex1, ex2}. In \cite{kh}, Khan showed some operational representations of $q$-polynomials. Koelink \cite{ko} explained some basic Lommel polynomials, Ismail \cite{is} considered the basic Bessel functions and polynomials. In \cite{rik}, Riyasat and Khan evidenced a determinant approach to $q$-Bessel polynomials and applications. In \cite{save}, Sahai and Verma established recursion formulas for $q$-Hypergeometric and $q$-Appell series.  In \cite{sa1, sa2}, Salem initiated and studied the theory of $q$-special matrix functions by introducing the matrix version of $q$-beta, $q$-gamma and $q$-hypergeometric matrix functions and its $q$-difference matrix equation. In \cite{sa3, sa4}, Salem considered the discrete basic (or $q$)-Hermite and $q$-Laguerre matrix polynomials. \cite{sh1, sh2} Shehata discussed the basic Konhauser and generalized Bessel matrix polynomials. In \cite{sr} Srivastava researched an introductory overview of Bessel, generalized Bessel and $q$-Bessel polynomials. In \cite{vesa} Verma and Sahai determined some recursion formulas for $q$-Lauricella series.
 
Motivated by these ideas, the $q$--BMPs and their basic (or $q$-)extensions have been studied extensively and widely. The objective of this study is to consider a class of $q$--BMPs from the viewpoint of their associated $q$-differential matrix equation, We also discuss the $q$-integral of these matrix polynomials. The product formula for the two $q$-Bessel matrix polynomials are obtained. Certain $q$-Laplace and $q$-Mellin transforms are investigated for the q-BMPs. Furthermore, the connections between $q$-Horn's matrix functions and $q$-Bessel matrix polynomials are discussed.
\subsection{Preliminaries}
In this subsection, we summarize an overview and an examination of the fundamental facts, notation, definitions, lemmas, and results of basic (or $q$-)special matrix functions that are needed for our work. Let $\mathbb{N}$ and $\Bbb{C}$ symbolize the sets of natural and complex numbers, and let $\Im(\xi)$ and  $\Re(\xi)$ symbolize the imaginary and real  parts of complex number $\xi$, and $0<|q|<1$, $q\in\Bbb{C}$. Let $\Bbb{C}^{\ell\times \ell}$ stand the vector space that contains all square matrices with $\ell$ columns and $\ell$ rows, where the elements are complex numbers. The $\mathbf{0}$ and $\mathbf{I}$ denote for the null and identity matrix in $\Bbb{C}^{\ell\times \ell}$.

For a matrix  $A$ in $\Bbb{C}^{\ell\times \ell}$, its spectrum $\sigma(A)$ denotes the set of all eigenvalues of $A$. We will use $||A||_{2}$ to the two-norm, which is defined by (see~\cite{j3, j4})
\begin{equation*}
\begin{split}
||A||_{2}=\sup_{x\neq 0}\frac{||A\chi||_{2}}{||\chi||_{2}},
\end{split}
\end{equation*}
where for a vector $\chi$ in $\Bbb{C}^{\ell}$, $||\chi||_2=(\chi^T\chi)^\frac{1}{2}$ is the Euclidean norm of $\chi$.
\begin{thm} If $\Omega(\xi)$ and $\Phi(\xi)$ are holomorphic functions of complex variable $\xi$, which are defined in an open set $\Theta$ of complex plane, then (see \cite{du})
\begin{equation*}
\begin{split}
\Phi(\mathbf{F})\Omega(\mathbf{E})=\Omega(\mathbf{E})\Phi(\mathbf{F}),
\end{split}
\end{equation*}
where $\mathbf{F}$, $\mathbf{E}$ are commutative matrices in $\Bbb{C}^{\ell\times \ell}$ with $\sigma(\mathbf{E}) \subset \Theta$ and $\sigma(\mathbf{F})\subset\Theta$.
\end{thm}
\begin{defn} For $\mathbf{E}$ in $\Bbb{C}^{\ell\times \ell}$, we say that $\mathbf{E}$ is a positive stable matrix if (see \cite{j3, j4})
\begin{equation}
\begin{split}
\Re(\omega) > 0,\quad \forall\; \omega\in \sigma(\mathbf{E}).\label{1.1}
\end{split}
\end{equation}
\end{defn}
The $q$-analogue of matrix $[\mathbf{E}]_{q}$ is defined as follows (see \cite{sa1}):
\begin{eqnarray}
\begin{split}
[\mathbf{E}]_{q}=\frac{I-q^{\mathbf{E}}}{1-q},q^{\mathbf{E}}=e^{\mathbf{E}\log q}, 0<|q|<1, q\in\mathbb{C}.\label{1.2}
\end{split}
\end{eqnarray}
The $q$-shifted factorial matrix function $(\mathbf{E};q)_{\upsilon}$ are defined as(see \cite{sa1})
\begin{eqnarray}
\begin{split}
(\mathbf{E};q)_{0}=I, (\mathbf{E};q)_{\upsilon}=\prod_{\tau=0}^{\upsilon-1}(I-q^{\tau}\mathbf{E}),||q^{\upsilon}\mathbf{E}||_{2}<1,\upsilon=1,2,\ldots,\label{1.3}
\end{split}
\end{eqnarray}
and
\begin{eqnarray}
\begin{split}
(\mathbf{E};q)_{\infty}=\lim_{\upsilon \rightarrow \infty}\prod_{\tau=0}^{\upsilon-1}(I-q^{\tau}\mathbf{E})=\prod_{\upsilon=0}^{\infty}(I-q^{\upsilon}\mathbf{E}),||q^{\upsilon}\mathbf{E}||_{2}<1,\label{1.4}
\end{split}
\end{eqnarray}
converges. Furthermore, if $\|q^{\upsilon}\mathbf{E}\|_{2}<1$ and $q^{-\upsilon}\notin\sigma(\mathbf{E}), \upsilon=0,1,2,\ldots$, the infinite product $(\mathbf{E};q)_{\infty}$ (\ref{1.4}) converges invertible matrix.
\begin{defn}
For any complex square matrix $\mathbf{F}$, the $q$-exponential matrix functions are defined as (see \cite{sa1})
\begin{eqnarray}
\begin{split}
E^{\mathbf{F}}_{q}=\sum_{\upsilon=0}^{\infty}\frac{q^{\binom{\upsilon}{2}}\mathbf{F}^{\upsilon}}{[\upsilon]_{q}!}=(-(1-q)\mathbf{F},q)_{\infty}\label{1.5}
\end{split}
\end{eqnarray}
and
\begin{eqnarray}
\begin{split}
e^{\mathbf{F}}_{q}=\sum_{n=0}^{\infty}\frac{\mathbf{F}^{n}}{[n]_{q}!}=[((1-q)\mathbf{F},q)_{\infty}]^{-1},\|\mathbf{F}\|_{2}<\frac{1}{|1-q|},(1-q)^{-1}q^{-\upsilon}\notin\sigma(\mathbf{F}), \upsilon=0,1,2,\ldots.\label{1.6}
\end{split}
\end{eqnarray}
\end{defn}
\begin{defn}
For any positive stable matrix $\mathbf{F}$, the $q$-Gamma matrix function $\Gamma_{q}(\mathbf{F})$ is defined as (see \cite{sa1}):
\begin{eqnarray}
\begin{split}
\Gamma_{q}(\mathbf{F})=\int_{0}^{\frac{1}{1-q}}t^{\mathbf{F}-I}E_{q}^{-qt}d_{q}t\label{1.7}
\end{split}
\end{eqnarray}
or
\begin{eqnarray}
\begin{split}
\Gamma_{q}(\mathbf{F})=\lim_{\upsilon \rightarrow \infty}(q,q)_{\upsilon}(1-q)^{I-\mathbf{F}}[(q^{\mathbf{F}},q)_{\upsilon}]^{-1}, q^{-\upsilon}\notin\sigma(q^{\mathbf{F}}),\upsilon=0,1,2,\ldots,.\label{1.8}
\end{split}
\end{eqnarray}
\end{defn}
\begin{defn}
For any positive stable matrix $\mathbf{E}$, the $q$-shifted factorial matrix function $(q^{\mathbf{E}},q)_{\upsilon}$ (\ref{1.3}) is defined as (see \cite{sa2})
\begin{eqnarray}
\begin{split}
(q^{\mathbf{E}},q)_{\upsilon}&=[\mathbf{E}]_{q}[\mathbf{E}+I]_{q}[\mathbf{E}+2I]_{q}\ldots[\mathbf{E}+(\upsilon-1)I]_{q}\\
&=(1-q)^{\upsilon}\Gamma_{q}{(\mathbf{E}+\upsilon I)}\Gamma^{-1}_{q}(\mathbf{E});\;
n\geq 1;\; (q^{\mathbf{E}},q)_{0}=I,\label{1.9}
\end{split}
\end{eqnarray}
\end{defn}
\begin{defn}
The $q$-Derivative operator $D_{\xi,q}$ of a function $\Psi$ at $\xi\neq0 \in \mathbb{C}$ is defined as (see \cite{ex2, ka})
\begin{equation}
D_{\xi,q}\Psi(\xi)= \frac{\Psi(\xi)-\Psi(q\xi)}{(1-q)\xi}, \xi\neq0,\label{1.11}
\end{equation}
and $D_{\xi,q}\Psi(0)=\frac{d\Psi(\xi)}{d\xi}|_{\xi=0}=f'(0)$, provided that $\Psi$ is differentiable at $\xi=0$, and~defined differential operator $\theta_{\xi,q}=\xi D_{\xi,q}$.
\end{defn}
\begin{defn}
Let $A$, $B$ be two complex square matrices and positive stable matrices, then a $q$-Beta matrix function $B_{q}(A,B)=$ is defined as follows (see \cite{sa1}):
\begin{eqnarray}
\begin{split}
B_{q}(A,B)=\int_{0}^{1}(tq,q)_{\infty}[(tq^{B},q)_{\infty}]^{-1}t^{A-I}d_{q}t,\label{1.12}
\end{split}
\end{eqnarray}
where $q^{-k}\notin\sigma(q^{A})\cup\sigma(q^{B})$ for all integer $k=0,1,2,\ldots,.$
\end{defn}
The $q$-Beta matrix function related to $q$-Gamma matrix function as follows: (see \cite{sa1})
\begin{eqnarray}
\begin{split}
B_{q}(\mathbf{F},\mathbf{E})=\Gamma_{q}{(\mathbf{E})}\Gamma_{q}(\mathbf{F})\Gamma^{-1}_{q}(\mathbf{F}+\mathbf{E}),\label{1.13}
\end{split}
\end{eqnarray}
where $\mathbf{F}$, $\mathbf{E}$ are commutative matrices, 
\begin{defn}
The $q$-Hypergeometric matrix functions $_{2}\phi_{1}$ is defined as (see \cite{sa2})
\begin{eqnarray}
\begin{split}
\;_{2}\phi_{1}\left(\mathbf{A},\mathbf{F};\mathbf{E};q,z\right) =
\sum_{\upsilon=0}^\infty\frac{z^{\upsilon}}{(q;q)_{\upsilon}} (\mathbf{A};q)_{\upsilon}(\mathbf{F};q)_{\upsilon}\big{[}(\mathbf{E};q)_{\upsilon}\big{]}^{-1},\label{1.14}
\end{split}
\end{eqnarray}
where $\mathbf{A}$, $\mathbf{F}$ and $\mathbf{E}$ are positive stable and commutative matrices in $\Bbb{C}^{N\times N}$ such that $q^{-\upsilon}\notin \sigma(\mathbf{E})$ for all integer $\upsilon\geq0$.
\end{defn}
\begin{lemma}
If $\Psi(\ell,\upsilon)$ is a matrix in $\Bbb{C}^{\ell\times \ell}$  for $\upsilon\geq 0$ and $\ell\geq 0$, then the relation is satisfied (see \cite{dw1}):
\begin{equation}
\begin{split}
\sum_{\upsilon=0}^{\infty}\sum_{\ell=0}^{\infty}\Psi(\ell,\upsilon)=\sum_{\upsilon=0}^{\infty}\sum_{\ell=0}^{\upsilon}\Psi(\ell,\upsilon-\ell).\label{1.15}
\end{split}
\end{equation}
According to (\ref{1.15}), we can write
\begin{equation}
\begin{split}
\sum_{\upsilon=0}^{\infty}\sum_{\ell=0}^{\upsilon}\Psi(\ell,\upsilon)=\sum_{\upsilon=0}^{\infty}\sum_{\ell=0}^{\infty}\Psi(\ell,\upsilon+\ell).\label{1.16}
\end{split}
\end{equation}
\end{lemma}
\begin{defn}
The $q$-Laplace transform of $f(x)$ is defined as (see \cite{ab,puk})
\begin{eqnarray}
\begin{split}
L_{s,q}\bigg{[}f(t)\bigg{]}=\frac{(q;q)_{\infty}}{s}\sum_{r=0}^{\infty}\frac{q^{r}}{(q;q)_{r}}f\bigg{(}\frac{q^{r}}{s}\bigg{)}.\label{1.17}
\end{split}
\end{eqnarray}
\end{defn}
\begin{defn}
The $q$-Mellin transform of $f(x)$ is defined as (see \cite{bo, uc})
\begin{eqnarray}
\begin{split}
M_{q}\bigg{[}f(x);s\bigg{]}&=(1-q)\sum_{r=0}^{\infty}q^{sr}f\big{(}q^{r}\big{)}.\label{1.18}
\end{split}
\end{eqnarray}
\end{defn}
\begin{defn}
For any finite positive integers $r$ and $s$, the generalized basic hypergeometric function is defined by (see \cite{dwsn})
\begin{equation}
\begin{split}
&\;_{r}\phi_{s}(q^{\mathbf{A}_{1}},q^{\mathbf{A}_{2}},\ldots,q^{\mathbf{A}_{r}};q^{\mathbf{B}_{1}},q^{\mathbf{B}_{2}},\ldots,q^{\mathbf{B}_{s}};z)\\
=&\sum_{k=0}^{\infty}\frac{z^{k}}{(q;q)_{\ell}}\prod_{\tau=1}^{r}(q^{\mathbf{A}_{\tau}};q)_{k}\bigg{[}\prod_{\upsilon=1}^{s}(q^{\mathbf{B}_{\upsilon}};q)_{k}\bigg{]}^{-1}\bigg{[}(-1)^{k}q^{\binom{k}{2}}\bigg{]}^{1+s-r},\label{1.19}
\end{split}
\end{equation}
where $\mathbf{A}_{\tau}$, $1 \leq \tau \leq r$ and $\mathbf{B}_{\upsilon}$, $1 \leq \upsilon \leq s$, are matrices in $\Bbb{C}^{\ell\times \ell}$ such that $q^{-k}\notin  \sigma(q^{\mathbf{B}_{\upsilon}})$, $1\leq \upsilon \leq s$ for all integers $k\geq 0$.
\end{defn}
\section{Derivation of main results}
\begin{defn}
Let $n$ be a positive integer and $A$ be a positive stable matrix in $\Bbb{C}^{\ell\times \ell}$, then we define the $q$-Bessel matrix polynomials $J_{n,q}(z;A)$ as
\begin{equation}
\begin{split}
J_{n,q}(z;A)=\frac{(q^{A-I};q)_{n}}{(q;q)_{n}}\;_{2}\Phi_{1}\bigg{(}q^{-nI},q^{A+(n-1)I};\mathbf{0};q,z\bigg{)}\label{2.1}
\end{split}
\end{equation}
\end{defn}
\begin{thm} 
The following the $q$-differential $q$-Hypergeometric matrix function $\;_{2}\Phi_{1}$ of second-order holds true :
\begin{eqnarray}
\begin{split}
\bigg{[}q^{A}q^{B}qz^{2}D^{2}_{q}+\bigg{(}z(q^{A}q^{B}+q^{A}[B]_{q}+q^{B}[A]_{q})-\frac{1}{1-q}I\bigg{)}D_{q}+[A]_{q}[B]_{q}\bigg{]}\;_{2}\Phi_{1}=\mathbf{0}\label{2.2}
\end{split}
\end{eqnarray}
\end{thm}
\begin{proof}
Applying the operator $\theta=zD=z\frac{d}{dz}$ for $q$-Hypergeometric matrix function $\;_{2}\Phi_{1}$, we get
\begin{eqnarray}
\begin{split}
[\theta I+A]_{q}\;_{2}\Phi_{1}\bigg{(}q^{A},q^{B};\mathbf{0};q,z\bigg{)}=[A]_{q}\;_{2}\Phi_{1}\bigg{(}q^{A+I},q^{B};\mathbf{0};q,z\bigg{)},\\
[\theta I+B]_{q}\;_{2}\Phi_{1}\bigg{(}q^{A},q^{B};\mathbf{0};q,z\bigg{)}=[B]_{q}\;_{2}\Phi_{1}\bigg{(}q^{A},q^{B+I};\mathbf{0};q,z\bigg{)},\\
[\theta]_{q}\;_{2}\Phi_{1}\bigg{(}q^{A},q^{B};\mathbf{0};q,z\bigg{)}=[A]_{q}[B]_{q}z\;_{2}\Phi_{1}\bigg{(}q^{A+I},q^{B+I};\mathbf{0};q,z\bigg{)}.\label{2.3}
\end{split}
\end{eqnarray}
From (\ref{2.3}), we have
\begin{eqnarray}
\begin{split}
\frac{1}{1-q}[\theta]_{q}I-z[\theta I+A]_{q}[\theta I+B]_{q}\;_{2}\Phi_{1}=0.\label{2.4}
\end{split}
\end{eqnarray}
Using the relation
\begin{eqnarray}
\begin{split}
[\theta I+A]_{q}&=[\theta]_{q}I+q^{\theta}[A]_{q}=[A]_{q}+q^{A}[\theta]_{q}I,\\
[\theta I+B]_{q}&=[\theta]_{q}I+q^{\theta}[B]_{q}=[B]_{q}+q^{B}[\theta]_{q} I,\label{2.5}
\end{split}
\end{eqnarray}
(\ref{2.4}) and simplifying, we obtain (\ref{2.2}).
\end{proof}
\begin{thm} The following $q$-differential $q$-Bessel matrix polynomials of second-order holds true:
\begin{eqnarray}
\begin{split}
\bigg{[}q^{A}z^{2}D^{2}_{q}+\bigg{(}\frac{z(q^{-nI}+q^{A+(n-1)I}-q^{A-I}-q^{A})-I}{1-q}\bigg{)}D_{q}+\frac{(I-q^{-nI})(I-q^{A+(n-1)I})}{(1-q)^{2}}\bigg{]}J_{n,q}(z;A)=\mathbf{0}.\label{2.6}
\end{split}
\end{eqnarray}
\end{thm}
\begin{proof}
Putting $A=-nI$ and $B=A+(n-1)I$ in (\ref{2.2}) and simplifying, we obtain (\ref{2.6}).
\end{proof}
\begin{thm} The following differential formula holds:
\begin{equation}
\begin{split}
D_{z,q}^{r}J_{n,q}(z;A)&=\frac{(q^{A-I};q)_{2r}}{(q-1)^{r}}q^{\binom{r}{2}-nr}J_{n-r,q}(z;A+2rI).\label{2.7}
\end{split}
\end{equation}
\end{thm}
\begin{proof}
Differentiating the equation (\ref{2.1}) with respect to $z$, yield
\begin{equation}
\begin{split}
D_{z,q}J_{n,q}(z;A)&=\frac{(q^{A-I};q)_{n}}{(q;q)_{n}}\frac{(I-q^{-nI})(I-q^{A+(n-1)I})}{1-q}\;_{2}\Phi_{1}\bigg{(}q^{-nI+I},q^{A+nI};\mathbf{0};q,z\bigg{)}\\
&=\frac{(q^{A-I};q)_{n}}{(q;q)_{n}}\frac{(I-q^{-nI})(I-q^{A+(n-1)I})}{1-q}\;_{2}\Phi_{1}\bigg{(}q^{-(n-1)I},q^{A+2I+(n-1)I-I};\mathbf{0};q,z\bigg{)}\\
&=\frac{(q^{A-I};q)_{n}}{(q;q)_{n}}\frac{(I-q^{-nI})(I-q^{A+(n-1)I})}{1-q}(q;q)_{n-1}[(q^{A+I};q)_{n-1}]^{-1}J_{n-1,q}(z;A+2I)\\
&=\frac{(I-q^{A-I})(I-q^{A})}{(q-1)}q^{-n}J_{n-1,q}(z;A+2I)=\frac{(q^{A-I};q)_{2}}{q-1}q^{-n}J_{n-1,q}(z;A+2I).\label{2.8}
\end{split}
\end{equation}
By repeating the above differentiation process $r$ times, we get
\begin{equation*}
\begin{split}
D_{z,q}^{r}J_{n,q}(z;A)&=\frac{(q^{A-I};q)_{2r}}{(q-1)^{r}}q^{-\{n+(n-1)+(n-2)+\cdots+(n-r+1)\}}J_{n-r,q}(z;A+2rI),
\end{split}
\end{equation*}
and using the relation
\begin{equation*}
\begin{split}
q^{-\{n+(n-1)+(n-2)+\cdots+(n-r+1)\}}=q^{-\{nr-(1+2+3+\cdots+(r-1))\}}=q^{-\{nr-\frac{r(r-1)}{2}\}}=q^{\binom{r}{2}-nr},
\end{split}
\end{equation*}
on simplifying, the desired result in (\ref{2.7}) is obtained.
\end{proof}
\begin{cor} Each of the following relations for $\;_{2}\Phi_{1}$ hold
\begin{eqnarray}
\begin{split}
\;_{2}\Phi_{1}\bigg{(}q^{A},q^{B};\mathbf{0};q,z\bigg{)}-\;_{2}\Phi_{1}\bigg{(}q^{A},q^{B};\mathbf{0};q,qz\bigg{)}&=(I-q^{A})(I-q^{B})z\;_{2}\Phi_{1}\bigg{(}q^{A+I},q^{B+I};\mathbf{0};q,z\bigg{)},\label{2.9}
\end{split}
\end{eqnarray}
\begin{eqnarray}
\begin{split}
\;_{2}\Phi_{1}\bigg{(}q^{A+I},q^{B};\mathbf{0};q,z\bigg{)}-\;_{2}\Phi_{1}\bigg{(}q^{A},q^{B};\mathbf{0};q,z\bigg{)}=zq^{A}(I-q^{B})\;_{2}\Phi_{1}\bigg{(}q^{A+I},q^{B+I};\mathbf{0};q,z\bigg{)},\label{2.10}
\end{split}
\end{eqnarray}
\begin{eqnarray}
\begin{split}
\;_{2}\Phi_{1}\bigg{(}q^{A},q^{B+I};\mathbf{0};q,z\bigg{)}-\;_{2}\Phi_{1}\bigg{(}q^{A},q^{B};\mathbf{0};q,z\bigg{)}=zq^{B}(I-q^{A})\;_{2}\Phi_{1}\bigg{(}q^{A+I},q^{B+I};\mathbf{0};q,z\bigg{)},\label{2.11}
\end{split}
\end{eqnarray}
\begin{eqnarray}
\begin{split}
\;_{2}\Phi_{1}\bigg{(}q^{A+I},q^{B};\mathbf{0};q,z\bigg{)}-\;_{2}\Phi_{1}\bigg{(}q^{A},q^{B+I};\mathbf{0};q,z\bigg{)}=z(q^{A}-q^{B})\;_{2}\Phi_{1}\bigg{(}q^{A+I},q^{B+I};\mathbf{0};q,z\bigg{)}\label{2.12}
\end{split}
\end{eqnarray}
and
\begin{eqnarray}
\begin{split}
q^{B}(I-q^{A})\;_{2}\Phi_{1}\bigg{(}q^{A+I},q^{B};\mathbf{0};q,z\bigg{)}-q^{A}(I-q^{B})\;_{2}\Phi_{1}\bigg{(}q^{A},q^{B+I};\mathbf{0};q,z\bigg{)}=(q^{B}-q^{A})\;_{2}\Phi_{1}\bigg{(}q^{A},q^{B};\mathbf{0};q,z\bigg{)}.\label{2.13}
\end{split}
\end{eqnarray}
\end{cor}
\begin{proof}
From the definition of the $\;_{2}\Phi_{1}$ and substitution on the L.H.S. of the equation (\ref{2.9}), we get
\begin{eqnarray*}
\begin{split}
&\;_{2}\Phi_{1}\bigg{(}q^{A},q^{B};\mathbf{0};q,z\bigg{)}-\;_{2}\Phi_{1}\bigg{(}q^{A},q^{B};\mathbf{0};q,qz\bigg{)}\\
&=\sum_{r=0}^{\infty}\frac{(1-q^{r})(q^{A};q)_{r}(q^{B};q)_{r}}{(q;q)_{r}}z^{r}=\sum_{r=1}^{\infty}\frac{(q^{A};q)_{r}(q^{B};q)_{r}}{(q;q)_{r-1}}z^{r}\\
&=\sum_{r=0}^{\infty}\frac{(q^{A};q)_{r+1}(q^{B};q)_{r+1}}{(q;q)_{r}}z^{r+1}=(I-q^{A})(I-q^{B})z\sum_{r=0}^{\infty}\frac{(q^{A+I};q)_{r}(q^{B+I};q)_{r}}{(q;q)_{r}}z^{r}\\
&=(I-q^{A})(I-q^{B})z\;_{2}\Phi_{1}\bigg{(}q^{A+I},q^{B+I};\mathbf{0};q,z\bigg{)}.
\end{split}
\end{eqnarray*}
Similar to the previous steps above, we have
\begin{eqnarray*}
\begin{split}
&\;_{2}\Phi_{1}\bigg{(}q^{A+I},q^{B};\mathbf{0};q,z\bigg{)}-\;_{2}\Phi_{1}\bigg{(}q^{A},q^{B};\mathbf{0};q,z\bigg{)}=\sum_{r=0}^{\infty}\bigg{(}(q^{A+I};q)_{r}-(q^{A};q)_{r}\bigg{)}\frac{(q^{B};q)_{r}}{(q;q)_{r}}z^{r}\\
&=\sum_{r=0}^{\infty}q^{A}(I-q^{rI})(I-q^{A})^{-1}\frac{(q^{A};q)_{r}(q^{B};q)_{r}}{(q;q)_{r}}z^{r}=\sum_{r=1}^{\infty}q^{A}(I-q^{A})^{-1}\frac{(q^{A};q)_{r}(q^{B};q)_{r}}{(q;q)_{r-1}}z^{r}\\
&=\sum_{r=0}^{\infty}q^{A}(I-q^{A})^{-1}\frac{(q^{A};q)_{r+1}(q^{B};q)_{r+1}}{(q;q)_{r}}z^{r+1}=q^{A}(I-q^{B})z\sum_{r=0}^{\infty}\frac{(q^{A+I};q)_{r}(q^{B+I};q)_{r}}{(q;q)_{r}}z^{r}\\
&=q^{A}(I-q^{B})z\;_{2}\Phi_{1}\bigg{(}q^{A+I},q^{B+I};\mathbf{0};q,z\bigg{)}.
\end{split}
\end{eqnarray*}
By repeating the previous steps above, we get (\ref{2.11}). From (\ref{2.10}) and (\ref{2.11}), we obtain the results (\ref{2.12}) and (\ref{2.13}).
\end{proof}
\begin{thm}
For the $q$-BMPs $J_{n,q}(z;A)$, we have the relations
\begin{equation}
\begin{split}
J_{n,q}(z;A)-J_{n,q}(qz;A)&=-(I-q^{A-I})(I-q^{A})q^{-n}zJ_{n-1,q}(z;A+2I),\label{2.14}
\end{split}
\end{equation}
\begin{eqnarray}
\begin{split}
(I-q^{A-I})J_{n-1,q}(z;A+I)-(1-q^{n})J_{n,q}(z;A)=zq^{-nI}(I-q^{A-I})(I-q^{A})J_{n-1,q}(z;A+2I),\label{2.15}
\end{split}
\end{eqnarray}
\begin{eqnarray}
\begin{split}
&(I-q^{A+(n-2)I})(I-q^{A-2I})^{-1}J_{n,q}(z;A-I)-J_{n,q}(z;A)\\
&=zq^{A+(n-1)I}(I-q^{A-I})(I-q^{A})(I-q^{A+(n-1)I})^{-1}J_{n-1,q}(z;A+2I),\label{2.16}
\end{split}
\end{eqnarray}
\begin{eqnarray}
\begin{split}
&(I-q^{A+(n-2)I})(I-q^{A+(n-1)I})\bigg{[}(I-q^{A+(n-3)I})J_{n-1,q}(z;A-I)-(1-q^{n})J_{n,q}(z;A-I)\bigg{]}\\
&=z(q^{-nI}-q^{A+(n-1)I})(I-q^{A-2I})(I-q^{A-I})(I-q^{A})J_{n,q}(z;A+2I)\label{2.17}
\end{split}
\end{eqnarray}
and
\begin{eqnarray}
\begin{split}
q^{A+(n-1)I}J_{n-1,q}(z;A+I)-q^{-nI}J_{n,q}(z;A+I)=(q^{A+(n-1)I}-q^{-nI})(I-q^{A-I})^{-1}J_{n,q}(z;A).\label{2.18}
\end{split}
\end{eqnarray}
\end{thm}
\begin{proof}
Substituting $A=-nI$, $B=A+(n-1)I$ in (\ref{2.9})-(\ref{2.13}) and simplifying, gives the desired results in (\ref{2.14})-(\ref{2.18}).
\end{proof}
\begin{thm}
The $q$-Bessel matrix polynomials $J_{n,q}(z;A)$ satisfy the pure recurrence relation as follows:
\begin{eqnarray}
\begin{split}
&\bigg{[}q^{A+(3n+1)I}+q^{(n+1)I}+q^{2(n+1)I}-q^{A+4nI}\bigg{]}J_{n+1,q}(z;A)=\\
&\bigg{[}q^{2A+(5n-1)I}+2q^{A+2nI}+2q^{A+(2n-1)I}-q^{2A+(4n-1)I}-q^{2A+4nI}-q^{2A+(3n-2)I}-q^{(n+1)I}\bigg{]}J_{n,q}(z;A)\\
&+\bigg{[}q^{2A+(4n-3)I}+2q^{2A+(4n-2)I}+q^{2A+(4n-1)I}+q^{2A+(3n-2)I}+q^{2A+(3n-1)I}+q^{2A+3nI}+q^{A+(n-2)I}\\
&+q^{A+(n-1)I}-q^{3A+(4n-3)I}-q^{3A+3(2n-1)I}-q^{2A+(2n-3)I}-2q^{A+2nI}-q^{A+(2n-1)I}-I\bigg{]}\;zJ_{n,q}(z;A)\\
&-\bigg{[}q^{3A+(5n-3)I}+q^{3A+4(n-1)I}+q^{A+(2n-1)I}-3q^{2A+(3n-2)I}-2q^{2A+3(n-1)I}-q^{2A+(4n-1)I}\bigg{]}J_{n-1,q}(z;A).\label{2.19}
\end{split}
\end{eqnarray}
\end{thm}
\begin{proof}
From the definition (\ref{2.1}), we can rewrite
\begin{equation}
\begin{split}
J_{n,q}(z;A)=\frac{(q^{A-I};q)_{n}}{(q;q)_{n}}\sum_{r=0}^{\infty}\frac{(q^{-nI};q)_{r}(q^{A+(n-1)I};q)_{r}}{(q;q)_{r}}z^{r} =\sum_{r=0}^{\infty}\Psi_{r}(A,n,q,z), \label{2.20}
\end{split}
\end{equation}
where $\Psi_{r}(A,n,q,z)=\frac{(q^{A-I};q)_{n}}{(q;q)_{n}}\frac{(q^{-nI};q)_{r}(q^{A+(n-1)I};q)_{r}}{(q;q)_{r}}z^{r}$.

In the same manner, we find that
\begin{equation}
\begin{split}
zJ_{n,q}(z;A)&=\frac{(q^{A-I};q)_{n}}{(q;q)_{n}}\sum_{r=0}^{\infty}\frac{(q^{-nI};q)_{r}(q^{A+(n-1)I};q)_{r}}{(q;q)_{r}}z^{r+1}\\
&=\frac{(q^{A-I};q)_{n}}{(q;q)_{n}}\sum_{r=1}^{\infty}\frac{(q^{-nI};q)_{r-1}(q^{A+(n-1)I};q)_{r-1}}{(q;q)_{r-1}}z^{r}\\
& =\sum_{r=0}^{\infty}(1-q^{r})(I-q^{-(n-r+1)I})^{-1}(I-q^{A+(n+r-2)I})^{-1}\Psi_{r}(A,n,q,z),\label{2.21}
\end{split}
\end{equation}
\begin{equation}
\begin{split}
J_{n+1,q}(z;A)&=\frac{(q^{A-I};q)_{n+1}}{(q;q)_{n+1}}\sum_{r=0}^{\infty}\frac{(q^{-(n+1)I};q)_{r}(q^{A+nI};q)_{r}}{(q;q)_{r}}x^{r}\\
&=-\sum_{r=0}^{\infty}q^{-(n+1)}(I-q^{-(n-r+1)I})^{-1}(I-q^{A+(n+r-1)I})\Psi_{r}(A,n,q,z)\label{2.22}
\end{split}
\end{equation}
and
\begin{equation}
\begin{split}
J_{n-1,q}(z;A)&=\frac{(q^{A-I};q)_{n-1}}{(q;q)_{n-1}}\sum_{r=0}^{\infty}\frac{(q^{-(n-1)I};q)_{r}(q^{A+(n-2)I};q)_{r}}{(q;q)_{r}}x^{r}\\
&=-q^{n}\sum_{r=0}^{\infty}(I-q^{-(n-r)I})(I-q^{A+(n+r-2)I})^{-1}\Psi_{r}(A,n,q,z).\label{2.23}
\end{split}
\end{equation}
For $n\geq1$, it follows that there exists a linear matrix recurrence relation in the form
\begin{eqnarray}
\begin{split}
J_{n+1,q}(z;A)=C_{1}J_{n,q}(z;A)+C_{2}\;zJ_{n,q}(z;A)-C_{3}J_{n-1,q}(z;A), \label{2.24}
\end{split}
\end{eqnarray}
where $C_{1}$, $C_{2}$ and $C_{2}$ are rational functions in $n$ and are independent of $z$.

Using equations (\ref{2.1}) and (\ref{2.20})-(\ref{2.23}) and reducing all coefficients of $\Psi_{r}(A,n,q,z)$ to the least common denominator, we see that
(\ref{2.24}) is equivalent to the following identity in $r$:
\begin{eqnarray*}
\begin{split}
-\sum_{r=0}^{\infty}q^{-(n+1)}(I-q^{-(n-r+1)I})^{-1}(I-q^{A+(n+r-1)I})\Psi_{r}(A,n,q,z)=C_{1}\bigg{(}\sum_{r=0}^{\infty}\Psi_{r}(A,n,q,z)\bigg{)}\\
+C_{2}\bigg{(}\sum_{r=0}^{\infty}(1-q^{r})(I-q^{-(n-r+1)I})^{-1}(I-q^{A+(n+r-2)I})^{-1}\Psi_{r}(A,n,q,z)\bigg{)}\\
-C_{3}\bigg{(}-q^{n}\sum_{r=0}^{\infty}(I-q^{-(n-r)I})(I-q^{A+(n+r-2)I})^{-1}\Psi_{r}(A,n,q,z)\bigg{)}.
\end{split}
\end{eqnarray*}
In the same manner, it follows that
\begin{eqnarray}
\begin{split}
-q^{-(n+1)}(I-q^{-(n-r+1)I})^{-1}(I-q^{A+(n+r-1)I})=C_{1}\\
+C_{2}\bigg{(}(1-q^{r})(I-q^{-(n-r+1)I})^{-1}(I-q^{A+(n+r-2)I})^{-1}\bigg{)}\\
-C_{3}\bigg{(}-q^{n}(I-q^{-(n-r)I})(I-q^{A+(n+r-2)I})^{-1}\bigg{)}.
\end{split}
\end{eqnarray}
Multiplying the whole equation by $ (I-q^{-(n-r+1)I})(I-q^{A+(n-r-2)I})$
\begin{eqnarray}
\begin{split}
-q^{-(n+1)}(I-q^{A+(n+r-1)I})(I-q^{A+(n+r-2)I})=C_{1}(I-q^{-(n-r+1)I})(I-q^{A+(n+r-2)I})\\
+C_{2}\bigg{(}(1-q^{r})\bigg{)}-C_{3}\bigg{(}-q^{n}(I-q^{-(n-r)I})(I-q^{-(n-r+1)I})\bigg{)}.
\end{split}
\end{eqnarray}
We need three equations for the determination of $C_{1}$, $C_{2}$ and $C_{3}$.\\
Coefficient of $q^{r}$, we get
\begin{eqnarray}
\begin{split}
q^{AI-3I}(1+q)=-C_{1}(q^{A+n-2}-q^{-n-1})-C_{2}-C_{3}(q+1).\label{2.27}
\end{split}
\end{eqnarray}
Coefficient of $(q^{r})^{2}$, we get
\begin{eqnarray}
\begin{split}
-q^{2A+n-4}=C_{1}(q^{A-3})+C_{3}(q^{-n-1}).\label{2.28}
\end{split}
\end{eqnarray}
Coefficient of $(q^{0})$, we get
\begin{eqnarray}
\begin{split}
-q^{-n-1}(I-q^{A+n-2}-q^{A+n-1}+q^{2A+2n-3})=C_{1}+C_{2}+C_{3}(q^{n}).\label{2.29}
\end{split}
\end{eqnarray}
By solving the three equations (\ref{2.27}), (\ref{2.28}) and (\ref{2.29}), we obtain the values of the constants $C_{1}$, $C_{2}$ and $C_{3}$.

Substituting the values thus obtained and clearing the result of fractions in (\ref{2.24}). Thus, the $q$-BMPs satisfy the pure matrix $q$-recurrence relation, we obtain (\ref{2.19}).
\end{proof}
\begin{thm}
The $q$-Laplace transform representation of $q$-Bessel matrix polynomials is given by
\begin{eqnarray}
\begin{split}
L_{s,q}\bigg{[}J_{n,q}(t;A)\bigg{]}=\frac{(q^{A-I};q)_{n}}{s(q;q)_{n}}\;_{3}\Phi_{2}\bigg{(}q^{-nI},q^{A+(n-1)I},q;\mathbf{0},\mathbf{0};q,\frac{1}{s}\bigg{)}\label{2.30}
\end{split}
\end{eqnarray}
and
\begin{eqnarray}
\begin{split}
L_{s,q}\bigg{[}t^{m-1}J_{n,q}(t;A)\bigg{]}=\frac{(q;q)_{\infty}(q^{A-I};q)_{n}}{s^{m}({q;q})_{n}(q^{m};q)_{\infty}}\;\;_{3}\Phi_{2}\bigg{(}q^{-nI},q^{A+(n-1)I},q^{mI};\mathbf{0},\mathbf{0};q,\frac{1}{s}\bigg{)},\label{2.31}
\end{split}
\end{eqnarray}
\end{thm}
\begin{proof}
By using the relations
\begin{eqnarray*}
\begin{split}
_{0}\Phi_{0}(-;q;t)=;\sum_{i=0}^{\infty}\frac{1}{(q;q)_i} t^{i}=\frac{1}{(t;q)_{\infty}},
\end{split}
\end{eqnarray*}
and
\begin{eqnarray*}
\begin{split}
_{0}\Phi_{0}(-;q;q^{m+k})=\frac{1}{(q^{m+k};q)_{\infty}}.
\end{split}
\end{eqnarray*}
By applying the $q$-Laplace transform formula (\ref{1.17}) of $q$-Bessel matrix polynomials and using the above relations, we get
\begin{eqnarray*}
\begin{split}
L_{s,q}\bigg{[}J_{n,q}(t;A)&\bigg{]}=\frac{(q;q)_{\infty}}{s}\sum_{r=0}^{\infty}\frac{q^{r}}{(q;q)_{r}}J_{n,q}(\frac{q^{r}}{s};A)\\
&=\frac{(q;q)_{\infty}}{s}\sum_{r=0}^{\infty}\frac{q^{r}}{(q;q)_{r}}\frac{(q^{A-I};q)_{n}}{(q;q)_{n}}\;_{2}\Phi_{1}\bigg{(}q^{-nI},q^{A+(n-1)I};\mathbf{0};q,\frac{q^{r}}{s}\bigg{)}\\
&=\frac{(q;q)_{\infty}}{s}\frac{(q^{A-I};q)_{n}}{(q;q)_{n}}\sum_{r=0}^{\infty}\sum_{s=0}^{\infty}\frac{q^{r}(q^{-nI};q)_{s}(q^{A+(n-1)I})_{s}}{(q;q)_{r}(q;q)_{s}}\bigg{(}\frac{q^{r}}{s}\bigg{)}^{s}\\
&=\frac{(q;q)_{\infty}}{s}\frac{(q^{A-I};q)_{n}}{(q;q)_{n}}\sum_{s=0}^{\infty}\frac{(q^{-nI};q)_{s}(q^{A+(n-1)I})_{s}}{(q;q)_{s}}\bigg{(}\frac{1}{s}\bigg{)}^{s}\sum_{r=0}^{\infty}\frac{q^{rs+r}}{(q;q)_{r}}\\
&=\frac{(q;q)_{\infty}}{s}\frac{(q^{A-I};q)_{n}}{(q;q)_{n}}\sum_{s=0}^{\infty}\frac{(q^{-nI};q)_{s}(q^{A+(n-1)I})_{s}}{(q;q)_{s}}\bigg{(}\frac{1}{s}\bigg{)}^{s}\;_{0}\Phi_{0}\bigg{(}-;\mathbf{0};q,q^{s+1}\bigg{)}\\
&=\frac{(q;q)_{\infty}}{s}\frac{(q^{A-I};q)_{n}}{(q;q)_{n}}\sum_{s=0}^{\infty}\frac{(q^{-nI};q)_{s}(q^{A+(n-1)I})_{s}}{(q;q)_{s}}\bigg{(}\frac{1}{s}\bigg{)}^{s}\frac{1}{(q^{s+1};q)_{\infty}}
\end{split}
\end{eqnarray*}
\begin{eqnarray*}
\begin{split}
&=\frac{(q;q)_{\infty}}{s}\frac{(q^{A-I};q)_{n}}{(q;q)_{n}}\sum_{s=0}^{\infty}\frac{(q^{-nI};q)_{s}(q^{A+(n-1)I})_{s}}{(q;q)_{s}}\bigg{(}\frac{1}{s}\bigg{)}^{s}\frac{(q;q)_{s}}{(q;q)_{\infty}}\\
&=\frac{1}{s}\frac{(q^{A-I};q)_{n}}{(q;q)_{n}}\sum_{s=0}^{\infty}\frac{(q^{-nI};q)_{s}(q^{A+(n-1)I})_{s}(q;q)_{s}}{(q;q)_{s}}\bigg{(}\frac{1}{s}\bigg{)}^{s}\\
&=\frac{1}{s}\frac{(q^{A-I};q)_{n}}{(q;q)_{n}}\;_{3}\Phi_{2}\bigg{(}q^{-nI},q^{A+(n-1)I},q;\mathbf{0},\mathbf{0};q,\frac{1}{s}\bigg{)}.
\end{split}
\end{eqnarray*}
By following the same steps above, we obtain (\ref{2.31}).
\end{proof}
\begin{thm} The $q$-Mellin transform representation of $q$-Bessel matrix polynomials is given by
\begin{eqnarray}
\begin{split}
M_{q}\bigg{[}J_{n,q}(z;A);s\bigg{]}=\frac{(q^{A-I};q)_{n}}{(q^{s};q)_{\infty}(q;q)_{n-1}}\;\;_{3}\Phi_{2}\bigg{(}q^{-nI},q^{A+(n-1)I},q^{sI};\mathbf{0},\mathbf{0};q,1\bigg{)}.\label{2.32}
\end{split}
\end{eqnarray}
\end{thm}
\begin{proof}
Using the definition of $q$-Mellin transform (\ref{1.18}) and applying of $q$-BMPs, we get
\begin{eqnarray*}
\begin{split}
M_{q}\bigg{[}J_{n,q}(z;A);s\bigg{]}&=(1-q)\sum_{r=0}^{\infty}q^{sr}J_{n,q}(q^{r};A)\\
&=(1-q)\frac{(q^{A-I};q)_{n}}{(q;q)_{n}}\sum_{r=0}^{\infty}q^{sr}\;\sum_{k=0}^{\infty}\frac{(q^{-nI};q)_k(q^{A+(n-1)I};q)_k}{(q;q)_k}q^{rk}\\
&=(1-q)\frac{(q^{A-I};q)_{n}}{(q;q)_{n}}\;\sum_{k=0}^{\infty}\frac{(q^{-nI};q)_k(q^{A+(n-1)I};q)_k}{(q;q)_k}\sum_{r=0}^{\infty}q^{sr+rk}\\
&=(1-q)\frac{(q^{A-I};q)_{n}}{(q;q)_{n}}\;\sum_{k=0}^{\infty}\frac{(q^{-nI};q)_k(q^{A+(n-1)I};q)_k}{(q;q)_k}\;_{0}\Phi_{0}\bigg{(}-;-;q,q^{s+k}\bigg{)}\\
&=(1-q)\frac{(q^{A-I};q)_{n}}{(q;q)_{n}}\;\sum_{k=0}^{\infty}\frac{(q^{-nI};q)_k(q^{A+(n-1)I};q)_k(q^{s};q)_{k}}{(q;q)_k(q^{s};q)_{\infty}}\\
&=\frac{(q^{A-I};q)_{n}}{(q^{s};q)_{\infty}(q;q)_{n-1}}\;\;_{3}\Phi_{2}\bigg{(}q^{-nI},q^{A+(n-1)I},q^{sI};\mathbf{0},\mathbf{0};q,1\bigg{)}
\end{split}
\end{eqnarray*}
\end{proof}
\begin{thm}
Each of the product formula for two  $q$-BMPs holds true
\begin{equation}
\begin{split}
&J_{n,q}(\alpha z;A)J_{m,q}(\beta z;A)=\frac{(q^{A-I};q)_{n}(q^{A-I};q)_{m}}{(q;q)_{n}(q;q)_{m}}\sum_{k=0}^{\infty}\sum_{\ell=0}^{k}\frac{(q^{-nI};q)_{k}(q^{A+(n-1)I};q)_{k}}{(q;q)_{k}}\alpha^{k}z^{k}\\
&\times\;_{3}\Phi_{2}\bigg{(}q^{-mI},q^{A+(m-1)I},q^{-k I};q^{(1-k+n)I},q^{(2-k-n)I-A};q,\frac{\beta}{\alpha}q^{(3-k) I-A}\bigg{)}.\label{2.33}
\end{split}
\end{equation}
\end{thm}
\begin{proof}
By using the formulas
\begin{eqnarray}
\begin{split}
(q^{-nI};q)_{k-\ell}&=(q^{-nI};q)_{k}[(q^{(1-k+n)I};q)_{\ell}]^{-1}(-q^{(n+1)I})^{\ell}q^{\binom{\ell}{2}-k\ell},\\
(q;q)_{k-\ell}&=\frac{(q;q)_{k}}{(q^{-k};q)_{\ell}}(-1)^{\ell}q^{\binom{\ell}{2}-k\ell},\\
(q^{A+(n-1)I};q)_{k-\ell}&=(q^{A+(n-1)I};q)_{k}[(q^{(2-k-n)I-A};q)_{\ell}]^{-1}(-q^{-A-(n-1)I}q)^{\ell}q^{\binom{\ell}{2}-k\ell}.\label{2.34}
\end{split}
\end{eqnarray}
Consider the product two $q$-Bessel matrix polynomials and using (\ref{1.15}) and (\ref{2.34}), becomes:
\begin{equation*}
\begin{split}
&J_{n,q}(\alpha z;A)J_{m,q}(\beta z;A)=\frac{(q^{A-I};q)_{n}(q^{A-I};q)_{m}}{(q;q)_{n}(q;q)_{m}}\sum_{k=0}^{\infty}\sum_{\ell=0}^{\infty}\frac{(q^{-nI};q)_{k}(q^{A+(n-1)I};q)_{k}}{(q;q)_{k}}\frac{(q^{-mI};q)_{\ell}(q^{A+(m-1)I};q)_{\ell}}{(q;q)_{\ell}}\alpha^{k}\beta^{\ell}z^{k+\ell}\\
&=\frac{(q^{A-I};q)_{n}(q^{A-I};q)_{m}}{(q;q)_{n}(q;q)_{m}}\sum_{k=0}^{\infty}\sum_{\ell=0}^{k}\frac{(q^{-nI};q)_{k-\ell}(q^{A+(n-1)I};q)_{k-\ell}}{(q;q)_{k-\ell}}\frac{(q^{-mI};q)_{\ell}(q^{A+(m-1)I};q)_{\ell}}{(q;q)_{\ell}}\alpha^{k-\ell}\beta^{\ell}z^{k}\\
&=\frac{(q^{A-I};q)_{n}(q^{A-I};q)_{m}}{(q;q)_{n}(q;q)_{m}}\sum_{k=0}^{\infty}\sum_{\ell=0}^{k}\frac{(q^{-nI};q)_{k}(q^{A+(n-1)I};q)_{k}}{(q;q)_{k}}(-1)^{\ell}q^{3\ell I-\ell A}\\
&\times q^{\binom{\ell}{2}-k\ell}\frac{(q^{-mI};q)_{\ell}(q^{A+(m-1)I};q)_{\ell}(q^{-k I};q)_{\ell}}{(q;q)_{\ell}}[(q^{(1-k+n)I};q)_{\ell}]^{-1}[(q^{(2-k-n)I-A};q)_{\ell}]^{-1}\alpha^{k-\ell}\beta^{\ell}z^{k}
\end{split}
\end{equation*}
\begin{equation*}
\begin{split}
&=\frac{(q^{A-I};q)_{n}(q^{A-I};q)_{m}}{(q;q)_{n}(q;q)_{m}}\sum_{k=0}^{\infty}\sum_{\ell=0}^{k}\frac{(q^{-nI};q)_{k}(q^{A+(n-1)I};q)_{k}}{(q;q)_{k}}\alpha^{k}z^{k}\frac{(q^{-mI};q)_{\ell}(q^{A+(m-1)I};q)_{\ell}(q^{-k I};q)_{\ell}}{(q;q)_{\ell}}\\
&\times[(q^{(1-k+n)I};q)_{\ell}]^{-1}[(q^{(2-k-n)I-A};q)_{\ell}]^{-1}\bigg{(}\frac{\beta}{\alpha}\bigg{)}^{\ell}(-1)^{\ell}q^{\binom{\ell}{2}}q^{((3-k) I-A)\ell}\\
&=\frac{(q^{A-I};q)_{n}(q^{A-I};q)_{m}}{(q;q)_{n}(q;q)_{m}}\sum_{k=0}^{\infty}\sum_{\ell=0}^{k}\frac{(q^{-nI};q)_{k}(q^{A+(n-1)I};q)_{k}}{(q;q)_{k}}\alpha^{k}z^{k}\\
&\times\;_{3}\Phi_{2}\bigg{(}q^{-mI},q^{A+(m-1)I},q^{-k I};q^{(1-k+n)I},q^{(2-k-n)I-A};q,\frac{\beta}{\alpha}q^{(3-k) I-A}\bigg{)}.
\end{split}
\end{equation*}
\end{proof}
\begin{thm} 
The integral representation for $q$-BMPs is given as
\begin{equation}
\begin{split}
J_{n,q}(z;A)=\frac{(1-q)^{n}\Gamma_{q}^{-1}(A-I)}{(q;q)_{n}}\int_{0}^{\frac{1}{1-q}}t^{A+(n-2)I}E_{q}^{-qt}\;_{1}\Phi_{0}\bigg{(}q^{-nI};-;q,(1-q^{n+1})tz\bigg{)}.\label{2.35}
\end{split}
\end{equation}
\end{thm} 
\begin{proof}
From (\ref{1.9}) and (\ref{2.1}), we have
\begin{equation*}
\begin{split}
&J_{n,q}(z;A)=\frac{(q^{A-I};q)_{n}}{(q;q)_{n}}\sum_{r=0}^{\infty}\frac{(q^{-nI};q)_{r}(q^{A+(n-1)I};q)_{r}}{(q;q)_{r}}z^{r}\\
&=\frac{1}{(q;q)_{n}}\sum_{r=0}^{\infty}\frac{(q^{-nI};q)_{r}(q^{A-I};q)_{n+r}}{(q;q)_{r}}z^{r}=\frac{\Gamma_{q}^{-1}(A-I)}{(q;q)_{n}}\sum_{r=0}^{\infty}\frac{(q^{-nI};q)_{r}}{(q;q)_{r}}\Gamma_{q}(A+(n+r-1)I)(1-q)^{n+r}z^{r}\\
&=\frac{(1-q)^{n}\Gamma_{q}^{-1}(A-I)}{(q;q)_{n}}\int_{0}^{\frac{1}{1-q}}t^{A+(n-2)I}E_{q}^{-qt}\sum_{r=0}^{\infty}\frac{(q^{-nI};q)_{r}}{(q;q)_{r}}(1-q^{n+1})^{r}(tz)^{r}d_{q}t\\
&=\frac{(1-q)^{n}\Gamma_{q}^{-1}(A-I)}{(q;q)_{n}}\int_{0}^{\frac{1}{1-q}}t^{A+(n-2)I}E_{q}^{-qt}\;_{1}\Phi_{0}\bigg{(}q^{-nI};-;q,(1-q^{n+1})tz\bigg{)}.
\end{split}
\end{equation*}
Thus the proof is completed.
\end{proof}
\begin{thm}
Each of the following integrals for $q$-Bessel matrix polynomial hold true
\begin{equation}
\begin{split}
&\frac{1}{1-q}\int_{0}^{1}t^{A-I}(1-qt)^{B-I}J_{n,q}(tz;A-I)d_{q}t=\frac{(q;q)_{\infty}(q^{A+B};q)_{\infty}(q^{A-2I};q)_{n}}{(q;q)_{n}}\\
&\times[(q^{A};q)_{\infty}]^{-1}[(q^{B};q)_{\infty}]^{-1}\;_{3}\phi_{2}(q^{-nI},q^{A},q^{A+(n-2)I};q^{A+B},\mathbf{0};q,z)\label{2.36}
\end{split}
\end{equation}
and
\begin{equation}
\begin{split}
&\frac{1}{1-q}\int_{0}^{1}t^{A-I}(1-qt)^{B-I}J_{n,q}((I-q^{B}t)z;A-I)d_{q}t=\frac{(q;q)_{\infty}(q^{A+B};q)_{\infty}(q^{A-2I};q)_{n}}{(q;q)_{n}}\\
&\times[(q^{A};q)_{\infty}]^{-1}[(q^{B};q)_{\infty}]^{-1}\;_{3}\phi_{2}(q^{-nI},q^{B},q^{A+(n-2)I};q^{A+B},\mathbf{0};q,z).\label{2.37}
\end{split}
\end{equation}
\end{thm}
\begin{proof}
From (\ref{2.1}) and (\ref{1.12}), we get
\begin{equation*}
\begin{split}
&\frac{1}{1-q}\int_{0}^{1}t^{A-I}(1-qt)^{B-I}J_{n,q}(tz;A-I)d_{q}t\\
&=\frac{(q^{A-2I};q)_{n}}{(1-q)(q;q)_{n}}\int_{0}^{1}t^{A-I}(1-qt)^{B-I}\sum_{r=0}^{\infty}\frac{(q^{-nI};q)_{r}(q^{A+(n-2)I};q)_{r}}{(q;q)_{r}}(tz)^{r}d_{q}t\\
&=\frac{(q^{A-2I};q)_{n}}{(1-q)(q;q)_{n}}\sum_{r=0}^{\infty}\frac{(q^{-nI};q)_{r}(q^{A+(n-2)I};q)_{r}}{(q;q)_{r}}z^{r}\int_{0}^{1}t^{A+(r-1)I}(1-qt)^{B-I}d_{q}t\\
&=\frac{(q^{A-2I};q)_{n}}{(q;q)_{n}}\sum_{r=0}^{\infty}\frac{(q^{-nI};q)_{r}(q^{A+(n-2)I};q)_{r}}{(q;q)_{r}}z^{r}(q;q)_{\infty}(q^{A+B+rI};q)_{\infty}[(q^{A+rI};q)_{\infty}]^{-1}[(q^{B};q)_{\infty}]^{-1}\\
&=\frac{(q;q)_{\infty}(q^{A+B};q)_{\infty}(q^{A-2I};q)_{n}}{(q;q)_{n}}[(q^{A};q)_{\infty}]^{-1}[(q^{B};q)_{\infty}]^{-1}\\
&\times\sum_{r=0}^{\infty}\frac{(q^{-nI};q)_{r}(q^{A};q)_{r}(q^{A+(n-2)I};q)_{r}}{(q;q)_{r}}[(q^{A+B};q)_{r}]^{-1}z^{r}.
\end{split}
\end{equation*}
Using the same techniques as proving equation ((\ref{2.36}), we obtain equation ((\ref{2.37}).
\end{proof}
\begin{thm} 
For the $q$-hypergeometric matrix functions $\;_{2}\Phi_{1}$, we have the following differentiation formulas
\begin{equation}
\begin{split}
D_{z,q}^{\mu}\bigg{[}z^{A+(\mu-1)I}\;_{2}\Phi_{1}\bigg{(}q^{A},q^{B};\mathbf{0};q,z\bigg{)}\bigg{]}=\frac{(q^{A};q)_{\mu}}{(1-q)^{\mu}}z^{A-I}\;_{2}\Phi_{1}\bigg{(}q^{A+\mu I},q^{B};\mathbf{0};q,z\bigg{)}\label{2.38}
\end{split}
\end{equation}
and
\begin{equation}
\begin{split}
D_{z,q}^{\mu}\bigg{[}z^{B+(\mu-1)I}\;_{2}\Phi_{1}\bigg{(}q^{A},q^{B};\mathbf{0};q,z\bigg{)}\bigg{]}=\frac{(q^{B};q)_{\mu}}{(1-q)^{\mu}}z^{B-I}\;_{2}\Phi_{1}\bigg{(}q^{A},q^{B+\mu I};\mathbf{0};q,z\bigg{)}.\label{2.39}
\end{split}
\end{equation}
\end{thm} 
\begin{proof}
Using the equation
\begin{equation*}
\begin{split}
\mathbf{D}_{q}^{\mu}z^{A+(\mu+r-1)I}=\frac{(q^{A};q)_{\mu}}{(1-q)^{\mu}}(q^{A+\mu I};q)_{r}[(q^{A};q)_{r}]^{-1}z^{A+(r-1)I},
\end{split}
\end{equation*}
on simplifying in L.H.S (\ref{2.38}), gives (\ref{2.38}). Similarly, we get (\ref{2.39}).
\end{proof}
\begin{thm}  For the $q$-BMPs $J_{n,q}(z;A)$, we have the following differentiation formulas
\begin{equation}
\begin{split}
D_{z,q}^{\mu}\bigg{[}z^{(\mu-n-1)I}J_{n,q}(z;A)\bigg{]}=\frac{(q^{A-I};q)_{\mu}}{(1-q)^{\mu}}q^{\binom{n}{2}-n\mu}z^{-(n+1)I}J_{n-\mu,q}(z;A+\mu I),z\neq0\label{2.40}
\end{split}
\end{equation}
and
\begin{equation}
\begin{split}
D_{z,q}^{\mu}\bigg{[}z^{A+(n+\mu-2)I}J_{n,q}(z;A)\bigg{]}=\frac{(q^{A-I};q)_{\mu}}{(1-q)^{\mu}}z^{A+(n-2)I}J_{n,q}(z;A+\mu I).\label{2.41}
\end{split}
\end{equation}
\end{thm} 
\begin{proof}
Substituting $A=-nI$, $B=A+(n-1)I$ in (\ref{2.38})-(\ref{2.39}) and using (\ref{2.1}), we obtain (\ref{2.40})-(\ref{2.41}).
\end{proof}
\begin{thm} For the $q$-hypergeometric matrix functions, we have the following differentiation relations
\begin{equation}
\begin{split}
D_{z,q}^{\mu}\bigg{[}z^{A+(\nu-1)I}\;_{2}\Phi_{1}\bigg{(}q^{A},q^{B};\mathbf{0};q,z\bigg{)}\bigg{]}=&\Gamma_{q}(A+\nu I)[\Gamma_{q}(A+(\nu-\mu)I)]^{-1}z^{A+(\nu-\mu-1)I}\\
&\times\;_{3}\Phi_{2}\bigg{(}q^{A+\mu I},q^{A+\nu I},q^{B};q^{A+(\nu-\mu)I},\mathbf{0};q,z\bigg{)}\label{2.42}
\end{split}
\end{equation}
and
\begin{equation}
\begin{split}
D_{z,q}^{\mu}\bigg{[}z^{B+(\nu-1)I}\;_{2}\Phi_{1}\bigg{(}q^{A},q^{B};\mathbf{0};q,z\bigg{)}\bigg{]}=&\Gamma_{q}(B+\nu I)[\Gamma_{q}(B+(\nu-\mu)I)]^{-1}z^{B+(\nu-\mu-1)I}\\
&\times\;_{3}\Phi_{2}\bigg{(}q^{A},q^{B+\mu I},q^{B+\nu I};q^{B+(\nu-\mu)I},\mathbf{0};q,z\bigg{)}.\label{2.43}
\end{split}
\end{equation}
\end{thm} 
\begin{proof}
Nota that
\begin{equation*}
\begin{split}
\mathbf{D}_{z,q}^{\mu}z^{A+(\nu+r-1)I}=\Gamma_{q}(A+\nu I)[\Gamma_{q}(A+(\nu-\mu)I)]^{-1}(q^{A+\nu I};q)_{r}[(q^{A+(\nu-\mu)I};q)_{r}]^{-1}z^{A+(\nu+r-\mu-1)I},
\end{split}
\end{equation*}
and
\begin{equation*}
\begin{split}
\Gamma_{q}(A+rI)=\frac{(q^{A};q)_{r}\Gamma_{q}(A)}{(1-q)^{r}}.
\end{split}
\end{equation*}
Using the above relations in L.H.S. (\ref{2.42})-(\ref{2.43}), we obtain (\ref{2.42})-(\ref{2.43}).
\end{proof}
\begin{thm} For the $q$-BMPs, we have the following differentiation relations
\begin{equation}
\begin{split}
D_{z,q}^{\mu}\bigg{[}z^{(\nu-n-1)I}J_{n,q}(z;A)\bigg{]}=&\frac{(q^{A-I};q)_{n}}{(q;q)_{n}}\Gamma_{q}((\nu-n) I)[\Gamma_{q}((\nu-\mu-n)I)]^{-1}z^{(\nu-\mu-n-1)I}\\
&\times\;_{3}\Phi_{2}\bigg{(}q^{(\mu-n) I},q^{(\nu-n) I},q^{A+(n-1)I};q^{(\nu-\mu-n)I},\mathbf{0};q,z\bigg{)}\label{2.44}
\end{split}
\end{equation}
and
\begin{equation}
\begin{split}
&D_{z,q}^{\mu}\bigg{[}z^{A+(n+\nu-2)I}J_{n,q}(z;A)\bigg{]}=\frac{(q^{A-I};q)_{n}}{(q;q)_{n}}\Gamma_{q}(A+(n+\nu-1)I)[\Gamma_{q}(A+(n+\nu-\mu-1)I)]^{-1}\\
&\times z^{A+(n+\nu-\mu-2)I}\;_{3}\Phi_{2}\bigg{(}q^{-nI},q^{A+(n+\mu-1)I},q^{A+(n+\nu-1)I};q^{A+(n+\nu-\mu-1)I},\mathbf{0};q,z\bigg{)}.\label{2.45}
\end{split}
\end{equation}
\end{thm} 
\begin{proof}
Substituting $A=-nI$, $B=A+(n-1)I$ in (\ref{2.42})-(\ref{2.43}) and using (\ref{2.1}), we get (\ref{2.44})-(\ref{2.45}).
\end{proof}
\begin{defn}
The $q$-Horn's matrix functions $H_{6}$ and $\Phi_{1}$ are defined by
\begin{eqnarray}
\begin{split}
&\Phi_{1}(q^{A},q^{B};q^{E};x,y)=\sum_{\tau,\upsilon=0}^{\infty}\frac{(q^{A};q)_{\tau+\upsilon}(q^{B};q)_{\tau}[(q^{E};q)_{\tau+\upsilon}]^{-1}}{(q;q)_{\tau}(q;q)_{\upsilon}}x^{\tau}y^{\upsilon};|x|<1,|y|<\infty\label{2.46}
\end{split}
\end{eqnarray}
and
\begin{eqnarray}
\begin{split}
&H_{6}(q^{A};q^{E};x,y)=\sum_{\tau,\upsilon=0}^{\infty}\frac{(q^{A};q)_{2\tau+\upsilon}[(q^{E};q)_{\tau+\upsilon}]^{-1}}{(q;q)_{\tau}(q;q)_{\upsilon}}x^{\tau}y^{\upsilon};|x|<\frac{1}{4},|y|<\infty,\label{2.47}
\end{split}
\end{eqnarray}
where $A$, $B$ and $E$ are matrices in $\Bbb{C}^{N\times N}$ satisfying the condition $q^{-(\upsilon+\tau)}\notin \sigma(E)$ for all integers $\upsilon+\tau\geq0$.
\end{defn}
\begin{thm} Let $A$, $B$ and $E$ be commutative matrices in $\Bbb{C}^{N\times N}$ satisfying the condition $q^{-\upsilon}\notin \sigma(E)$ for all integers $\upsilon\geq0$, then connections between the $q$-Horn's matrix functions $H_{6}$, $\Phi_{1}$ and $q$-hypergeometric matrix function satisfy the following 
\begin{eqnarray}
\begin{split}
H_{6}(q^{A};q^{E};-xy,y)=&\sum_{m=0}^{\infty}\frac{(-1)^{m}(q^{A};q)_{2m}[(q^{E};q)_{m}]^{-1}}{(q;q)_{m}}x^{m}y^{m}\\
&\times\;_{2}\Phi_{1}\bigg{(}q^{-mI},\mathbf{0};q^{(1-2m)I-A};q,\frac{q^{(1-m)I-A}}{x}\bigg{)}\label{2.48}
;|xy|<\frac{1}{4}.
\end{split}
\end{eqnarray}
and
\begin{eqnarray}
\begin{split}
\Phi_{1}(q^{A},q^{B};q^{E};-xy,y)=\sum_{n=0}^{\infty}\frac{(q^{A};q)_{n}[(q^{E};q)_{n}]^{-1}}{(q;q)_{n}}y^{n}\;_{2}\Phi_{0}\bigg{(}q^{-nI},q^{B};-;q,q^{n}x\bigg{)};|xy|<1.\label{2.49}
\end{split}
\end{eqnarray}
\end{thm}
\begin{proof}
Observe that, we have 
\begin{eqnarray}
\begin{split}
(q;q)_{m-n}&=\frac{(q;q)_{m}}{(q^{-m};q)_{n}}q^{\binom{n}{2}-mn},\;(q;q)_{n-m}=\frac{(q;q)_{n}}{(q^{-n};q)_{m}}q^{\binom{m}{2}-mn},\\
(q^{A};q)_{2m-n}&=(q^{A};q)_{2m}[(q^{(1-2m)I-A};q)_{n}]^{-1}(-q^{I-A})^{n}q^{\binom{n}{2}-2mn}\\
&=(q^{A};q)_{m}(q^{A+mI};q)_{m}[(q^{(1-2m)I-A};q)_{n}]^{-1}(-q^{I-A})^{n}q^{\binom{n}{2}-2mn},\\
(q^{B};q)_{m-n}&=(q^{B};q)_{m}[(q^{(1-m)I-B};q)_{n}]^{-1}(-q^{I-B})^{n}q^{\binom{n}{2}-mn}.\label{2.50}
\end{split}
\end{eqnarray}
From (\ref{2.47}), (\ref{2.50}), (\ref{1.15}) and (\ref{2.1}), we have
\begin{eqnarray*}
\begin{split}
&H_{6}(q^{A};q^{E};-xy,y)=\sum_{m=0}^{\infty}\sum_{n=0}^{m}\frac{(-1)^{m-n}(q^{A};q)_{2m-n}[(q^{E};q)_{m}]^{-1}}{(q;q)_{m-n}(q;q)_{n}}x^{m-n}y^{m}\\
&=\sum_{m=0}^{\infty}\sum_{n=0}^{m}\frac{(-1)^{m-n}[(q^{E};q)_{m}]^{-1}}{(q;q)_{n}}\frac{(q^{A};q)_{2m}}{(q^{(1-2m)I-A};q)_{n}}(-q^{I-A})^{n}q^{-mn}\frac{(q^{-m};q)_{n}}{(q;q)_{m}}x^{m-n}y^{m}\\
&=\sum_{m=0}^{\infty}\frac{(-1)^{m}(q^{A};q)_{2m}[(q^{E};q)_{m}]^{-1}}{(q;q)_{m}}x^{m}y^{m}\;_{2}\Phi_{1}\bigg{(}q^{-mI},\mathbf{0};q^{(1-2m)I-A};q,\frac{q^{(1-m)I-A}}{x}\bigg{)}.
\end{split}
\end{eqnarray*}
Using the same steps above, we can prove (\ref{2.49}).
\end{proof}
\begin{thm} The connections between $q$-Horn's matrix functions $H_{6}$, $\Phi_{1}$ and $q$-BMPs satisfy the following
\begin{eqnarray}
\begin{split}
H_{6}(q^{A};q^{E};-xy,y)=\sum_{n=0}^{\infty}[(q^{E};q)_{n}]^{-1}y^{n}q^{-\binom{n}{2}}J_{n;q}(x;A+I) ;|xy|<\frac{1}{4}\label{2.51}
\end{split}
\end{eqnarray}
and
\begin{eqnarray}
\begin{split}
\Phi_{1}(q^{A},q^{B};q^{E};-xy,y)&=\sum_{m=0}^{\infty}(-1)^{m}(q^{A};q)_{m}(q^{B};q)_{m}x^{m}[(q^{E};q)_{m}]^{-1}[(q^{(1-2m)I-B};q)_{n}]^{-1}\\
&\times J_{m;q}\bigg{(}2(1-m)I-B,\frac{q^{I-B}y}{x}\bigg{)};|xy|<1,\label{2.52}
\end{split}
\end{eqnarray}
where $A$, $B$ and $E$ are commutative matrices in $\Bbb{C}^{N\times N}$ satisfying the condition $q^{-\upsilon}\notin \sigma(E)$ for all integers $\upsilon\geq0$.
\end{thm}
\begin{proof} 
The L.H.S. of (\ref{2.51}) is equal to
\begin{eqnarray*}
\begin{split}
&H_{6}(q^{A};q^{E};-xy,y)=\sum_{n=0}^{\infty}\sum_{m=0}^{n}\frac{(-1)^{m}(q^{A};q)_{n+m}[(q^{E};q)_{n}]^{-1}}{(q;q)_{m}(q;q)_{n-m}}x^{m}y^{n}\\
&=\sum_{n=0}^{\infty}\frac{(q^{A};q)_{n}[(q^{E};q)_{n}]^{-1}}{(q;q)_{n}}y^{n}q^{-\binom{n}{2}}\sum_{m=0}^{n}\frac{(q^{-nI};q)_{m}}{(q;q)_{m}}(q^{A+nI};q)_{m}\frac{}{}(-xq^{n})^{m}\\
&=\sum_{n=0}^{\infty}[(q^{E};q)_{n}]^{-1}y^{n}q^{-\binom{n}{2}}J_{n;q}(x;A+I).
\end{split}
\end{eqnarray*}
which lead to (\ref{2.51}). In the same way as proving (\ref{2.51}), we obtain (\ref{2.52}).
\end{proof}
\section{Concluding remarks}
In the present investigation, the family constructed in this paper ensures us various substantial applications for $q$-Bessel matrix polynomials ($q$-BMPs). Certain properties comprehensive matrix $q$-recurrence relations, $q$-difference equations and $q$-integral representations, $q$-Laplace transforms and $q$-Mellin transforms of the $q$-BMPs have been established under certain conditions on matrices. The product formula for two $q$-BMPs are derived. The connections between the $q$-Horn’s matrix functions $H_{6}$, $\Phi_{1}$, $q$-hypergeometric matrix function and $q$-BMPs are discussed, which will form the basis for future work.

\end{document}